\input{psfig}
\documentstyle[12pt]{article}
\def\forces{\parallel\!\!\! -}

\def\hexnumber#1{\ifcase#1 0\or1\or2\or3\or4\or5\or6\or7\or8\or9\or
        A\or B\or C\or D\or E\or F\fi }

\font\teneuf=eufm10
\font\seveneuf=eufm7
\font\fiveeuf=eufm5
\newfam\euffam
\textfont\euffam=\teneuf
\scriptfont\euffam=\seveneuf
\scriptscriptfont\euffam=\fiveeuf
\def\frak{\fam\euffam \teneuf}

\font\tenmsx=msam10
\font\sevenmsx=msam7
\font\fivemsx=msam5
\font\tenmsy=msbm10
\font\sevenmsy=msbm7
\font\fivemsy=msbm5
\newfam\msxfam
\newfam\msyfam
\textfont\msxfam=\tenmsx  \scriptfont\msxfam=\sevenmsx
  \scriptscriptfont\msxfam=\fivemsx
\textfont\msyfam=\tenmsy  \scriptfont\msyfam=\sevenmsy
  \scriptscriptfont\msyfam=\fivemsy
\edef\msx{\hexnumber\msxfam}

\mathchardef\upharpoonright="0\msx16
\let\restriction=\upharpoonright
\def\Bbb#1{\tenmsy\fam\msyfam#1}

\def\qed{{\vcenter{\hrule height.4pt \hbox{\vrule width.4pt height5pt
 \kern5pt \vrule width.4pt} \hrule height.4pt}}}
\def\notin{{\in}\kern-5.5pt / \kern1pt}
\def\ok{\vbox{\hrule height 8pt width 8pt depth -7.4pt
    \hbox{\vrule width 0.6pt height 7.4pt \kern 7.4pt \vrule width 0.6pt height 7.4pt}
    \hrule height 0.6pt width 8pt}}
\def\nt{{\leq}\kern-1.5pt \vrule height 6.5pt width.8pt depth-0.5pt \kern 1pt}
\def\sd{{\times}\kern-2pt \vrule height 5pt width.6pt depth0pt \kern1pt}
\def\zp#1{{\hochss Y}\kern-3pt$_{#1}$\kern-1pt}

\def\BB{{\Bbb B}}

\def\LL{{\Bbb L}}

\def\PP{{\Bbb P}}
\def\QQ{{\Bbb Q}}
\def\RR{{\Bbb R}}

\def\11{{\Bbb 1}}

\font\capit=cmcsc10 scaled\magstep0

\newtheorem{theorem}{Theorem}[section]
\newtheorem{proposition}[theorem]{Proposition}
\newtheorem{corollary}[theorem]{Corollary}
\newtheorem{definition}{Definition}
\newtheorem{lemma}[theorem]{Lemma}
\newtheorem{problem}{Open Question}
\newtheorem{question}{Question}

\title{Heights of Models of $ZFC$ and the Existence of End Elementary
Extensions}
\author{Andr\'es Villaveces}

\begin{document}

\maketitle
\begin{abstract}
\noindent
The existence of End Elementary Extensions of models M
of ZFC is related to the ordinal height of M, according to classical
results due to Keisler, Morley and Silver. In this paper, we further
investigate the connection between the height of M and the existence of
End Elementary Extensions of M. In particular, we prove that the theory
`ZFC + GCH + there exist measurable cardinals + all
inaccessible non weakly compact cardinals are possible heights of models
with no End Elementary Extensions' is
consistent relative to the theory `ZFC + GCH + there exist measurable
cardinals + the weakly compact cardinals are cofinal in ON'.
We also provide a simpler coding that destroys GCH but otherwise
yields the same result.
\end{abstract}

\bigskip

\noindent
I wish to thank my advisor, Kenneth Kunen, for many helpful
conversations and comments, always full of interesting insights.
I also wish to thank Ali Enayat for very helpful discussions about some of the
topics treated in this article, Sy Friedman for interesting questions
and comments related to the class forcing construction around Theorem
\ref{theorem:heightandgch}, Mirna D\v zamonja and Arnie Miller
for various helpful discussions.
\bigskip

\noindent

\section{Introduction.}

\medskip
\noindent
Let $(M,E),(N,F)$, etc. denote models of `enough set theory.'
The central notion of extension we use in this article is the
well known `end elementary extension'.. A model $(N,F)$
{\bf end extends} $(M,E)$ iff for every $a\in M$, the sets
$a_E = \{ b\in M | b E a\}$ and $a_F = \{ b\in N | b F a\}$ are the same.
In other words, elements of $M$ are not enlarged by the extension
from $M$ to $N$.

\medskip
\noindent
Let $({\cal E}_M,\prec _e)$ denote the structure of all non-trivial
elementary end extensions (`eees') of $M$ (a model of set theory), together
with the relation $\prec _e$ (we write `${\cal A}\prec _e{\cal B}$' if and
only if ${\cal B}$ is an
elementary {\it end\/} extension of ${\cal A}$). $\prec _e$ is an ordering on
${\cal E}_M$.
The kind of ordering relation that $\prec _e$ is on ${\cal E}_M$ depends
heavily on certain structural features of $M$. Notice that
structures like $({\cal E}_M,\prec _e)$ need not be well-founded: as Kaufmann
notes in [Ka 83], if $\kappa$ is weakly compact, $({\cal
E}_{{\cal R}(\kappa )},\prec _e)$ has infinite descending chains.

\bigskip
\noindent
In this paper, we concentrate on problems of {\bf existence} of end
elementary extensions of specific models of ZFC, and we study the
relationship between the height of those models and the
possibility to obstruct the existence of eees. In a forthcoming work,
we will concentrate on the study of chains in $({\cal E}_M,\prec _e)$,
and their connection with large cardinal properties of $M$.

\bigskip
\noindent
We consider the following general question:

\begin{question}\label{general}
How does the structure of $M$ affect the structure of
$({\cal E}_M,\prec _e)$?
\end{question}

\medskip
\noindent
The earliest results toward a solution to Question~\ref{general} were
obtained by Keisler, Silver and Morley in [KeSi 70] and [KeMo 68]. All
their theorems addressed the specific question of the {\bf
existence} of eees (when is ${\cal E}_M \not= 0$?). They provided answers for
the two following cases.

\begin{theorem} (Keisler, Morley~[KeMo 68])
Let $M$ be a model of ZFC, $cof(M) = \omega$. Then ${\cal E}_M \not= 0$.
(This takes care in particular of all countable models of ZFC.)
\end{theorem}

\begin{theorem}\label{keislersilver} (Keisler, Silver~[KeSi 70])
If $M$ is of the form ${\cal R}(\kappa )$, where $\kappa$
is a weakly compact cardinal, then for all $S \subset M$, ${\cal
E}_{({\cal R}(\kappa ),\in ,S)} \not= 0$.
\end{theorem}

\noindent
In light of those early results, it is natural
to ask to what degree does the {\bf height} of $M$ determine
the structure of ${\cal E}_M$. We begin by the basic question: does the height
of $M$ determine the {\bf existence} of an eee?

\bigskip
\noindent
Here are the two main results of this paper.

\bigskip
\noindent
{\bf Theorem \ref{theorem:height}:}
 The theory `ZFC + $\exists \lambda$ measurable + $\forall \kappa (\kappa $
 inaccessible not weakly compact $\to \exists$ transitive $M_\kappa \models
ZFC$ such that $o(M)=\kappa$ and ${\cal E}^{wf}_M=0$)' is consistent relative
to the theory `ZFC + $\exists \lambda$ measurable'.

\bigskip
\noindent
{\bf Theorem \ref{theorem:heightandgch}:}
The theory `ZFC + GCH + $\exists \lambda (\lambda$ measurable) + $\forall
\kappa[\kappa$ inaccessible not weakly compact $\to  \exists$ transitive
$M_\kappa \models
ZFC$ such that $o(M)=\kappa$ and ${\cal E}^{wf}_M=0$)]' is
consistent relative to the theory `ZFC + $\exists \lambda (\lambda$
measurable) + the weakly compact cardinals are cofinal in ON'.

\bigskip
\noindent
The `conclusion' in both theorems is that the existence of
transitive models of the form ${\cal R}(\kappa)$
which have eees, yet the inner model $M_\kappa \subset {\cal R}(\kappa)$
does not have any well-founded eees (that is, the `global' negative answer
to the Height Problem, later abbreviated as the `NED' property) is
consistent with fairly reasonable large cardinal axioms.

\bigskip
\noindent
The proofs in both cases use codings of the obstructions to the existence of
eees {\bf by an appropriate inner model}. Although the results look similar,
the codings used in both cases differ strongly: in the first theorem, the
places where GCH holds or fails provide the main tool for the coding; in the
second case, since one needs to obtain models where GCH holds
everywhere, one cannot use anymore such a device for the coding. In that case,
forcing nonreflecting stationary sets at the appropriate levels to a model
previously freed of any of those sets does the trick.

\bigskip
\noindent
The notation we use is standard. We restrict ourselves to the
case of {\it transitive\/} $M$. Following two different traditions, we
freely switch between the two notations
`$V_\kappa$' and `${\cal R}(\kappa )$' when we denote the set of objects of
the universe of rank less than $\kappa$. Given a model $M$, $o(M)$ denotes
the ordinal height of $M$.

\bigskip

\section{The Height Problem.}

\medskip
\noindent
When $M$ is definable, the problem is trivial, by the following

\begin{proposition} If $M$ is definable in ${\cal R(\kappa)}$, and
${\cal R(\kappa )}$ has an eee, then $M$ has an eee.
\end{proposition}

\noindent
{\bf Proof:} Let $N$ elementarily end extend ${\cal R(\kappa)}$. Since $o(M) = 
\kappa$, $M^N$, the interpretation of the definition of $M$ in $N$, properly
(end) extends $M$. It is easy to verify that we actually have $M {\prec}_e
 M^N$. \hfill $\Box$
\medskip

\noindent
Of course, this is far from settling the general problem

\begin{question}\label{question:height}
{\bf (Height Problem)} Do ${\cal E}_{\cal R(\kappa)}\not= 0$ and $o(M) =
\kappa$ imply together that ${\cal E}_M\not= 0$?
\end{question}

\begin{figure}[h]\hspace*{50mm}
\psfig{figure=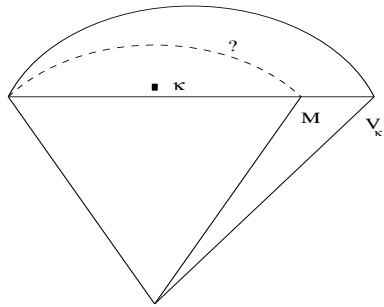,height=4cm,width=5cm}
\caption{ The Height Problem }    
\end{figure}

\noindent
As Theorem \ref{theorem:height} shows, it is quite possible that this question
have a negative answer: there are models of set theory where a certain
$\cal R(\kappa)$ has eees, yet some transitive submodel with the same
height has no eees. This can even be obtained in a quite homogeneous
way: the `extendability property' only holds at weakly compact cardinals.

\medskip
\noindent
Still, if we add new axioms to the theory,
the answer to the question may become positive. This is the case, for example,
when $V=L$
holds in $M$, the question becomes trivial, since in that case the only
transitive submodel of height $\kappa$ of an $\cal R(\kappa)$ is
$\cal R(\kappa)$ itself.

\medskip
\noindent
Therefore, the answer to the Height Problem is independent of ZFC.

\bigskip

\subsection{Consistency of a global negative answer.}

\noindent
We will next present some of the possible situations
for answers to a version of the problem. Let $G(\lambda )$ denote any large
cardinal property preserved under Easton-type extensions where the iteration
process is not carried too often (measurability, etc.). More precisely,
we mean here properties unaffected by Easton iterations which only act at
certain {\bf successor} cardinals.

\begin{definition}
Let $\kappa$ be an inaccessible cardinal. We say that
$\kappa$ is NON-END-DETERMINING (or for
short $NED(\kappa )$) if and only if there is a transitive $M_\kappa \models
ZFC$ such that $o(M)=\kappa$ and ${\cal E}^{wf}_M=0$.
\end{definition}

\noindent
Thus, if $NED(\kappa )$, there may be models of the form ${\cal R}_\kappa$
which have eees, yet certain inner models of them do not have
well-founded eees. Observe
that the $NED$ property depends strongly on the universe where $\kappa$ is
being considered.

\begin{theorem}\label{theorem:height}
 The theory `ZFC + $\exists \lambda G(\lambda )$ + $\forall \kappa (\kappa $
 inaccessible not weakly compact $\to$ NED($\kappa$))' is consistent relative
to the theory `ZFC + $\exists \lambda G(\lambda )$'.
\end{theorem}

\noindent
We will devote the remainder of this section to proving Theorem
\ref{theorem:height} and to draw some corollaries from it. 
The proof of Theorem \ref{theorem:height} is based on a forcing
construction which produces codes for
various sets which witness the non existence of eees at the desired
cardinals. The `$\exists
\lambda G(\lambda )$' clause is added to the theory to point out that the
forcing can be done in such a way that large cardinal properties
$G(\lambda )$ (at least, those not destroyed by Easton iterations
which only act at certain {\bf successor} cardinals) are consistent with the
negative answer to Question \ref{question:height}.

\medskip
\noindent
In the next lemma, we introduce the forcing construction, and we show that
both weak compactness and non weak compactness are preserved in the extension.

\begin{lemma}\label{pres} Let $M$ be a model of ZFC+GCH, and let $M[G]$ be a
forwards Easton extension obtained by adding $\alpha ^{++}$ subsets
to every successor inaccessible $\alpha$. Then $M$ and $M[G]$ have the
same weakly compact cardinals.
\end{lemma}

\noindent
{\bf Proof:} Let first $\kappa$ be weakly compact. Notice that $\kappa$ is not 
in the domain of the Easton iteration: it only acts at successor
inaccessibles. Let $\phi$ be
${\Sigma}_1^1$ (say $\phi \equiv \exists X\psi (X)$, $\psi$ first order). Let
$$1\forces _{\PP} (\forall \alpha <\kappa   \exists X\subset \alpha
 [(\alpha ,<,\tau )\models \psi (X)]),$$
for some $\PP$-name $\tau$ in $V$ such that $1\forces \tau \subset \kappa$.
Then, by the $\Pi _1^1$-indescribability of $\kappa$ in $M$,
there are $N$, $\PP^*$, $\tau ^*$ such that
$$({\cal R}(\kappa ),\in ,\tau ,\PP ){\prec }_e (N,\in ,\tau ^*, \PP ^*).$$
Now, in $M$, $\PP ^*$ is an Easton forcing which (since $\PP ^*$ does not
act on $\kappa$) can be seen as a product
$$\PP ^* = \PP  \times  \PP '.$$
A usual argument shows that $\PP '$ is $(< \kappa ^+)$-closed. Then, $\PP ^*$
does not add any new subsets to $\kappa$. Now, in $M$, we have
$$1\forces _{\PP ^*} \exists X \subset  \kappa  \Big( (\kappa ,<,\tau ,X)
\models  \psi \Big) .$$
So, since $M^{\PP ^*} = (M^{\PP '})^{\PP }$,
$$1\forces _{\PP '} \exists \mbox { some }\PP\mbox {-name }\sigma  \Big[ 1
\forces _{\PP}
 \Big( (\kappa ,<,\tau ,\sigma )\models  \psi \Big) \Big].$$
Hence, by the $\kappa$-closure of $\PP '$, a fusion argument in the
style of Silver
provides a $\PP$-name $\sigma$ in $M$ such that
$$1\forces _{\PP} \Big( (\kappa ,<,\tau ,\sigma )\models  \psi \Big).$$
So, $\kappa$ is also weakly compact in $M[G]$.

\noindent
Now, if $\kappa$ is {\it not\/} weakly compact in the universe, there is
a $\kappa$-Aronszajn tree $T$. There are two cases:

\begin{description}
\item[ (i) ] $\kappa$ is Mahlo: then $\PP$ is the product of a
$\kappa$-Knaster and a $\kappa ^+$-closed forcing; thus
$$1\forces \mbox { T is }\kappa \mbox {-Aronszajn,}$$
whereby $\kappa$ cannot be weakly compact in $M[G]$.

\item[ (ii) ] $\kappa$ is not Mahlo. Then, there is a club $C$ of singular
cardinals
in $\kappa$, $C \in M$. None of these singular cardinals becomes regular
in $M[G]$; clearly $C$ is also a club in $\kappa$ in $M[G]$. So, $\kappa$
remains non-Mahlo in the extension. \hfill $\Box_{\ref{pres}}$
\end{description}

\medskip
\noindent
Having shown this preservation property of $\PP$, we are ready to complete
the
\medskip

\noindent
{\bf Proof of Theorem \ref{theorem:height}:} Assume that GCH
and $\exists \lambda
(G(\lambda ))$ both hold in $V$. Let $V[G]$ be a forward Easton
extension like in the previous Lemma. Then,
$V$ and $V[G]$ have the same weakly compact cardinals. If $\kappa $ is
{\it not\/} weakly compact (in $V$ or $V[G]$), then we need to find
a transitive $M_{\kappa}$ in $V[G]$ with no eee and $o(M_\kappa )=\kappa$.
There are three cases:

\begin{description}
\item[ 1) ] $\kappa$ is a successor
inaccessible: then let $M_{\kappa} = {\cal R}(\kappa )$. This $M_{\kappa}$
cannot be end elementary extended: when we have $N\succ _eM_\kappa$, $\kappa$
turns out to be inaccessible in $N$. But this implies that the inaccessibles
are unbounded in $\kappa$, and this is incompatible with the fact that
$\kappa$ is a successor inaccessible.

\item[ 2) ] $\kappa$ is
a non-Mahlo limit inaccessible: then fix $C$ a club in $\kappa$ with no
regular cardinals, and let $M_{\kappa}$ code $C$ by the powers of
the successor
inaccessibles below $\kappa$. This $M_{\kappa}$ is obtained as follows: write
$\kappa$ as an increasing sequence $(\lambda _{\alpha})_{\alpha <\kappa}$,
where the $\alpha$'s run over the successor inaccessibles below $\kappa$.
Now, set
$M_{\kappa}$ = the extension of ${\cal R}(\kappa )$ where $\alpha ^{++}$
subsets of $\alpha$ are added only when $\lambda_{\alpha}\in C.$
Clearly, we have ${\cal R}(\kappa )^V
\subset M_{\kappa} \subset {\cal R}(\kappa )^{V[G]}$ and $o(M_\kappa )=
\kappa$. Since $M_{\kappa}$
codes $C$, it cannot have any eee (if it had one, say $N$, then $\kappa$
would be a singular cardinal in $N$!).

\item[ 3) ] $\kappa$ is
Mahlo: then let $T$ be a $\kappa$-Aronszajn tree in $V$. Since $T$ remains
$\kappa$-Aronszajn in $V[G]$, it is enough to let $M_{\kappa}$ encode
$T$ similarly to part 2; this new $M_{\kappa}$ cannot be elementarily end
extended,
for it would then provide a $\kappa$-branch to $T$ in $V[G]$, contradicting
the fact that $T$ is Aronszajn. This completes our proof. \hfill
$\Box _{\ref{theorem:height}}$
\end{description}

\medskip
\begin{corollary}
The Height Problem may have a negative answer.
\end{corollary}

\noindent
{\bf Proof:} In the model $V[G]$ constructed in Theorem \ref{theorem:height},
let $\kappa$ be the first inaccessible such that ${\cal E}_{\cal R(\kappa )}
\not= 0$. Results in [KeSi 70] show that $\kappa$ is not weakly
compact.
Then, by Theorem \ref{theorem:height}, there exists a transitive model
$M_\kappa$
of $ZFC$ which is of height $\kappa$ and has no well-founded eees.
This provides a negative
answer to the Height Problem. \hfill $\Box$

\medskip
\noindent
It is interesting to observe the following, in connection to
Theorem~\ref{theorem:height}.

\begin{proposition}
If $\kappa$
{\bf is} weakly compact, then for all transitive $M \subset
{\cal R}(\kappa )$ of height $\kappa$, ${\cal E}_M\not= 0$.
\end{proposition}

\noindent
{\bf Proof:} It suffices to observe that if $\kappa$ is
weakly compact, then $({\cal R}(\kappa ),\in ,S)$ has eees, for
all $S\subset {\cal R}(\kappa )$. Take any (transitive) $M \subset
{\cal R}(\kappa )$ of height $\kappa$. Then $({\cal R}(\kappa ),\in ,M)$
has a well founded eee $(N,\in ,M')$. Clearly,

\begin{description}
\item[ i ] $M\subset _eM'$: by endness, all the `new' elements go `on top'
of M, and
\item[ ii ] $M\prec M'$: given any sentence $\phi$, ${\cal R}(\kappa )
\models (M\models \phi)$ if and only if $N\models (M'\models \phi)$. This
implies $M\prec M'$ by the transitivity of $N$. \hfill $\Box$
\end{description}

\medskip
\noindent
The `opposite' problem is still open, for cardinals $\kappa$ which are
inaccessible not weakly compact:

\begin{problem}
Is the theory `ZFC + large cardinals + $\forall \kappa$ inaccessible
not weakly compact
$((M\models ZFC \wedge o(M) = \kappa \wedge {\cal E}_{\cal R(\kappa )} 
\not= 0) \to {\cal E}_M \not= 0)$' 
consistent?
\end{problem}

\noindent
In other words, is it consistent to have large cardinals and
simultaneously a {\it globally\/} positive answer to the Height Problem?
Up to now, the only way to get globally positive answers to the
Height Problem is with $V=L$ or similar axioms.

\medskip
\noindent
{\bf Remark:} An easy modification of the construction given for Theorem
\ref{theorem:height} provides models for the Height Problem only at certain
cardinals.

\medskip
\noindent
We basically have two extreme situations here: the Height Problem under
the very restrictive $V=L$ (trivial positive answer since there are
no strictly inner models), and the Height
Problem in the presence of Large Cardinals (negative answers).
It is natural to ask

\begin{question}
What happens `in between'\/ $V=L$ and Large Cardinal Axioms? More
specifically, what is the
situation for models of the form $L[0^\sharp ]$ or $L[\mu ]$ (where $\mu$
is a measure in the usual sense)?
\end{question}

\noindent
The bottom line here is that in the presence of large cardinals, there are
enough sharps which allow encoding the parameters on a model of the
form ${\cal R}(\kappa )$, for $\kappa$ non weakly compact. So, a situation
similar to that of Theorem \ref{theorem:height} would be obtained via 
a $\kappa$-closed poset: in $L[\mu ]$, if $\kappa$ is not weakly compact,
then there is a model $M$ of height $\kappa$ with no eees. Situations similar
to these are the subject of the final section.

\bigskip
\noindent
\section{A new coding: keeping GCH true throughout the forcing.}

\medskip
\noindent
In the previous section, the construction of models for $NED(\kappa )$
used Mc Aloon's coding method:
coding (in inner models of the generic extension) the Aronszajn trees or
the clubs of singulars needed to destroy eees by using the successor
inaccessibles where GCH
holds/fails as the coding device. Naturally, if one wants to obtain a similar
result while {\it keeping\/} GCH true, one must code the construction
in a completely different way.

\begin{theorem}\label{theorem:heightandgch}
The theory `ZFC + GCH + $\exists \lambda (\lambda$ measurable) + $\forall
\kappa[\kappa$ inaccessible not weakly compact $\to$ NED($\kappa)]$' is
consistent relative to the theory `ZFC + $\exists \lambda (\lambda$
measurable) + the weakly compact cardinals are cofinal in ON'.
\end{theorem}

\noindent
This answers a question that Sy Friedman asked me during the Tenth Latin
American Mathematical Logic Symposium in Bogot\'a.
The proof of this theorem consists
of a two step forcing, followed by a construction of models of height
$\kappa$ with no eees in the same way as in the end of the proof of Theorem
\ref{theorem:height}: coding in $M_\kappa$ an object that makes it
impossible for $M_\kappa$ to have eees. (Depending on the case, a
$\kappa$-Aronszajn tree or a club of singulars.)

\medskip
\noindent
The two forcing constructions provide enough coding tools to complete the
proof in a way analogous to how coding was used to prove Theorem
\ref{theorem:height}. The idea is to
use the existence/non-existence of nonreflecting stationary sets on
$\alpha^{++}$ (those which will be called
$E_{\alpha ^{++}}^\alpha$-sets in the construction), for $\alpha$ successor
inaccessible, instead of the continuum
function (cardinals where GCH holds/fails), to do the coding.

\bigskip

\begin{figure}[ht]\hspace*{8.00mm}\label{picture:proof}
\psfig{figure=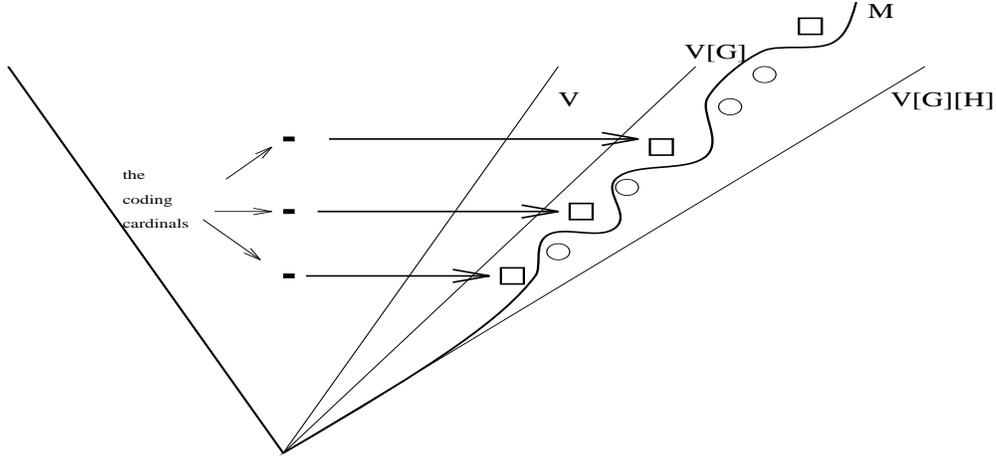,height=6cm,width=13cm}
\caption{ Coding M by adding only {\bf some} E-sets (the $\Box$s) }
\end{figure}

\medskip
\noindent
We start with
a model of `$ZFC + GCH + \exists \lambda (\lambda$ measurable) + the weakly
compacts are cofinal in $ON$'. We first force the
non-existence of `$E$-sets', and then force back the existence of them on
successor inaccessibles. We will then be able to construe the inner models
we need inside the final generic
model by putting in it {\bf exactly} the $E$-sets one needs to carry the
coding.

\begin{definition} ($E_\kappa ^\theta$, nonreflecting stationary sets)
Let $\theta$ be a regular cardinal, and $\kappa >\theta$, with $cof(\kappa )
>\aleph _1$. We mean by $E_\kappa ^\theta$ that
$$\exists A \subset \{\delta <\kappa | cof(\delta ) = \theta\},\mbox {
A stationary in }\kappa \mbox { such that }$$
$$\forall \mbox { limit }\alpha <\kappa (A\cap \alpha \mbox { is not
stationary in }\alpha).$$
Abusing the notation a bit, we sometimes say that the set $A$ is an `$E_\kappa
^\theta$-set', or even an `E-set'.
\end{definition}

\medskip
\noindent
{\bf Proof of Theorem \ref{theorem:heightandgch}:} Start with $M,$ a model
of `ZFC + GCH + $\exists \lambda (\lambda$ measurable) + $\forall \alpha
\exists \mu >\alpha (\mu$ weakly compact)'.

\noindent
Let $\PP$ be the Easton product of the L\'evy
collapses $\LL _i=Coll(\kappa_i ^+,<\theta_i )$ of cardinals $<\theta _i$ to
$\kappa ^+_i$, where $i\in ORD^M$, and
the sequence $\big<<\kappa _i,\theta _i>|i\in ORD^M \big>$ is
defined as follows:
 $\kappa_0$
is the first successor inaccessible cardinal, and $\theta _0$ is the first
weakly compact cardinal; if $\kappa _i$, $\theta _i$ are defined, then
$\kappa _{i+1}$ is the first successor inaccessible {\it after\/}
$\theta _i$ and $\theta _{i+1}$ is the first weakly compact cardinal after
$\kappa _{i+1}$. For limit $i$, $\kappa _i$ is the first successor
inaccessible after $\sup \{ \kappa _j | j<i\}$, and $\theta _i$ is the
first weakly compact cardinal after $\kappa _i$.

\noindent
More precisely, $p\in \PP$ iff $p$ is a function
such that
\begin{description}
\item[ i ] $p(i)\in \LL _i$,
\item[ ii ] for all $i\in$ dom($p$), for all $\sigma$ inaccessible such
that $i<\sigma$ implies $\theta_i<\sigma$,
dom($p$)$\cap \sigma$ is bounded in $\sigma$.
\end{description}
\noindent
$\PP$ is ordered in the usual way.

\noindent
For every $\kappa$, $\PP$ factors into three forcings: $\PP \approx
\PP_{<\kappa }\ast \PP_\kappa \ast
\PP_{>\kappa }$, where $\PP_{<\kappa }$, $\PP_\kappa$, $\PP_{>\kappa}$
mean respectively the portions of the iteration below $\kappa$, at $\kappa$,
and above $\kappa$.  
\noindent
By Easton's Lemma and $GCH$ in $M$, $\PP$ collapses the cardinals in
$[\kappa _i^+,\theta _i)$ to $\kappa _i^+$ (thus, in the generic extension,
$\theta _i$ becomes $\kappa _i^{++}$), and preserves all the other cardinals.
In particular, all the $\kappa _i$s remain successor inaccessible in the
extension, as the reader can easily check.

\medskip
\noindent
{\bf Claim:} For all $i$, for all stationary $S \subset \{ \beta
<\kappa _i^{++} |$ cof$(\beta )= \kappa _i\}$ in $V^\PP$, there exists $\sigma$
such that cof$(\sigma ) = \kappa _i^{+}$ and $S\cap \sigma$ is stationary.
(We say that the
stationary reflection property holds for $\kappa _i$ and $\kappa _i^{++}$ at
$\kappa _i^+$.)

\noindent
To show this, it is enough to look at $V^{(\PP _{\leq \theta _i})}$, since
this coincides with $V^\PP$ up to the level $\kappa _{i+1}$. $V^{(\PP _{\leq
\theta _i})}$ can be seen as $(V^{\LL _i})^{\PP_{<\theta _i}}$. Since the
stationary reflection property holds in $V$ for $\kappa _i$ and $\theta _i^+$
at $\kappa _i^+$ (by Baumgartner [Ba 76]), we have
that the stationary reflection property for $\kappa _i$ and $\theta _i =
\kappa _i^{++}$ (at $\kappa _i^+$ necessarily)
also holds in $V^{\LL _i}$. On the other hand, $|\PP _{<\theta _i}|
< \kappa _i^+ <\theta _i$; hence if $S$ is a stationary subset of $\{ \beta
<\kappa _i^{++} |$ cof$(\beta )= \kappa _i\}$, and $S\in (V^{\LL
_i})^{\PP_{<\theta _i}}$, then there is $T\subset S$ stationary,
$T\in V^{\LL _i}$. By the stationary reflection property in $V^{\LL _i}$,
there exists some $\gamma <\theta _i$, with cof$(\gamma) = \kappa _i^+$,
such that $T\cap \gamma$ is stationary in $V^{\LL _i}$. $T\cap \gamma$ is
still stationary in $(V^{\LL _i})^{\PP_{<\theta _i}}$, since $|\PP_{<\theta
_i}| < \kappa _i^+$. (Since whenever $p \forces `\dot \tau$ club in
$\gamma$' (cof $\gamma$ = $\kappa _i^+$), we have that $\dot \tau$ already
is a club in $\gamma$ in
$V^{\LL _i}$.) All this ends the proof of the Claim.

\noindent
So, the stationary reflection property holds in $V^\PP$ for all $\kappa _i$,
$\theta _i$. This ends the first stage of the proof. We have obtained a model
$V_1=V[G]$, for $G$ $\PP$-generic, where $E_{\theta _i} ^{\kappa _i}$-sets
do NOT exist, for any $i\in ORD$.

\bigskip
\noindent
For the second step of the proof, we will simplify matters by writing `$V$'
instead of `$V_1$' and calling our forcing $\PP$ again. We need to force
$E_{\kappa ^{++}} ^\kappa$-sets into
$\kappa ^{++}$, for each successor inaccessible $\kappa$. We will achieve
this goal by `piecewise adding' the $E$-sets.

\noindent
The proof that such $E$-sets can be added is an Easton iteration.
We present here a variant of the construction by Cummings, D\v zamonja and
Shelah in [CuD\v zSh 95] of a generic extension of a model of `$GCH + \exists
\lambda (\lambda$ measurable)', where a regular cardinal $\theta$ `strongly
non-reflects' at $\lambda _i$, where $i<\theta$, cof($\lambda _i$)=$\lambda
_i<\lambda _i^+ < \theta$. Their construction yields in particular $E_\theta
^{\lambda _i}$-sets in the generic extension (see Lemma 2 and Theorem 1
of [CuD\v zSh 95]).

\medskip
\noindent
{\bf The Building Blocks:} given any successor inaccessible cardinal $\kappa$,
let $\PP _\kappa$ consist of conditions $p=(\alpha _p,E_p,C_p)$, where
$\alpha _p<\kappa^{++}$, cof($\alpha _p$)=$\kappa ^+$, $E_p\subset\alpha _p$,
$$C_p: S_{\kappa ^+}^{\alpha _p}\to
\bigcup _{\beta <\alpha _p}{\cal P}(\beta)$$

\noindent
where $S_{\kappa ^+}^{\alpha _p} = \{ \beta <\alpha _p |\mbox
{ cof}(\beta )=\kappa ^+\}$, and for all $\beta \in S_{\kappa
^+}^{\alpha _p}$, $C_p(\beta)$ is a club on $\beta$. We also require that
for all $\beta \in S_{\kappa ^+}^{\alpha _p}$, $C_p(\beta )\cap E_p = 0$,
and $E_p\subset S_\kappa^{\alpha _p}$.
So, $E_p$ behaves like an `approximation' to a non-reflecting stationary
subset of $S_\kappa^{\kappa^{++}}$.
We define the ordering of $\PP _\kappa$ by

\bigskip
\qquad \qquad \qquad
$q\leq p \longleftrightarrow \left \{
\begin{tabular}{l}
(a) $\alpha _q \geq \alpha _p$,\\
(b) $E_q\cap \alpha _p=E_p$,\\
(c) $C_q\restriction S^{\kappa ^+}_{\alpha _p} = C_p$.
\end{tabular}\right .$

\bigskip

\noindent
We need only worry at ordinals of cofinality $\kappa ^+$:
\begin{lemma} If cof$(\alpha )\in [\aleph _1,\kappa^+)$, then
there is a club $C\subset \alpha$ such that for all $\gamma \in C$,
cof$(\gamma )<\kappa$.
\end{lemma}

\noindent
{\bf Proof:} Fix $\alpha$ with $\aleph _1\leq$ cof$(\alpha )\leq \kappa$, and
let $f:\mbox {cof}(\alpha) \to \alpha$ be increasing and
continuous.
Then all the $\gamma$'s in the domain of $f$ are $<\kappa$.
As $f$ is continuous increasing, ran($f$) is a club in $\alpha$, and
for every $\beta \in$ ran($f$), we have cof($\beta) =$ cof$(f^{-1}(\beta ))
<\kappa$. \hfill $\Box$

\medskip
\noindent
Let now $G$ be $\PP _\kappa$-generic, and let $E_G = \bigcup _{p\in  G}E_p$.

\medskip
\noindent
{\bf Claim:} $E_G$ is a nonreflecting stationary subset of $S^\kappa _{\kappa
^{++}}$.

\noindent
To show that $E_G$ is nonreflecting, let $\alpha <\kappa^{++}$ be of
cofinality $\kappa ^+$. Then $\alpha <\alpha _p$ for some $p\in G$, and
thus $C_{\alpha _p}(\alpha )\cap E_{\alpha _p} = 0$. But by (b) in the
definition of $\leq$, $E_G\cap\alpha_p
=E_{\alpha _p}$, $C_{\alpha _p}(\alpha )\subset \alpha <\alpha _p$. So,
$C_{\alpha _p}(\alpha )$, which is a
club on $\alpha$,
does not intersect $E_G$: $C_{\alpha _p}(\alpha 
)\cap E_G = (C_{\alpha _p}(\alpha )\cap \alpha _p)\cap E_G =
C_{\alpha _p}(\alpha )\cap (\alpha _p\cap E_G) = C_{\alpha _p}(\alpha )\cap
E_{\alpha _p} = 0$.

\noindent
To see that $E_G$ is
stationary, let $\Gamma$ be a canonical name for $E_G$, and suppose that
$p\forces (\tau \mbox { club of }\kappa^{++} \wedge \tau \cap \Gamma = 0)$. Get then a
decreasing $\kappa$-chain of conditions $p>
p_0>p_1>\dots >p_\xi >\dots$, $p_\xi \in \PP _\kappa$, such that for $\xi<\kappa$, $p_\xi \forces
\delta _\xi \in \tau$ and $\alpha _{p_\xi }<\delta _\xi <\alpha _{p_{\xi +1}}$.
Let $q=(\sup _{\xi <\kappa}\alpha _{p_\xi },\bigcup _{\xi
<\kappa}E_{p_\xi },\bigcup _{\xi <\kappa}C_{p_\xi })$. This $q$ is not a
condition, since the cofinality of its support is $<\kappa ^+$. But we
just have to take any condition $r$ that extends $q$ in a natural
way.

\noindent
[How? Just let
$r = \Big(\sup _{\xi <\kappa}\alpha _{p_\xi }+\kappa ^+,\bigcup _{\xi
<\kappa}E_{p_\xi }\cup \{ \alpha ^* \}
,\bigcup _{\xi <\kappa}C_{p_\xi }\cup \hat C\Big)$, where
\begin{itemize}
\item [ \bf i ] $\alpha ^* = \sup _{\xi <\kappa}\alpha _{p_\xi }$, and
\item[ \bf ii ] $\hat C(\beta )$ is any club on $\beta$ with
min($\hat C(\beta ))>\alpha ^*$, for $\beta \in
S^{\alpha ^*+\kappa ^+}_{\kappa ^+}\setminus
\alpha ^*$.
\end{itemize}

\noindent
It is not difficult to check that this $r\in \PP _\kappa$.]

\medskip
\noindent
Then, the limit of the $\delta _\xi$'s ($= \alpha ^*$) is forced by
$r$ to belong to
the club $\tau$. On the other hand, $\alpha ^*\in E_r$ and $r\forces
\alpha ^*\in \Gamma$. But this contradicts the fact that $p\forces (\tau
\mbox {
club }\wedge \tau \cap \Gamma = 0)$.

\medskip
\noindent
{\bf Claim:} $\PP _\kappa$ is $(<\kappa ^+)$-closed.

\noindent
This Claim is established by the same argument as the previous one,
although taking now sequences of arbitrary lengths below $\kappa ^+$.

\medskip
\noindent
{\bf Remark:} Clearly, $\PP _\kappa$ is not $(<\kappa ^{++})$-closed, since
$\PP _\kappa$ adds $\kappa ^{++}$-sequences.

\noindent
The $\PP _\kappa$'s correspond thus to the `Building Blocks' of our
construction. The remaining part consists of the iteration through all the
successor inaccessible cardinals, and of the proof that measurables are
preserved.

\bigskip
\noindent
{\bf Claim:} $|\PP _\kappa | \leq \kappa ^{++}$.

\noindent
(Just observe that in $V_1$,
$\kappa ^{++}\cdot 2^{\kappa ^+}\cdot 2^{\kappa ^+} = \kappa ^{++}$.)

\bigskip
\noindent
\subsection{The Iteration.}

\noindent
The iteration is in the style of Backward Easton. The
supports are bounded below
regular cardinals. We add new subsets only at successor inaccessible stages.
Formally, this corresponds to defining $\PP_{\leq \alpha}$ as the forcing
up to stage $\alpha$ and $\dot\QQ _\alpha \in V^{\PP _\alpha}$ as the
forcing at stage $\alpha$. Set $\dot\QQ _\alpha = \{0\}$, if $\alpha$
is not a successor inaccessible, and $\dot\QQ _\alpha = (\PP
_\alpha)^{V_{<\alpha}}$, otherwise, where $V_{<\alpha}$ stands for
$V^{\PP _{<\alpha }}$.
By the $\kappa ^+$-closure of the forcing $\PP _\kappa$, this iteration
yields a
model of $ZFC$. In order to finish the proof, we need only show that the
measurability of $\lambda$ is preserved by this iteration.

\medskip
\noindent
Let $\lambda$ be a measurable cardinal in $V$, and let $j:V\to M$ be the
ultrapower map arising from a normal measure $U$ on $\kappa$. We first
observe that it is enough to prove that $\lambda$ remains measurable in
the extension by $\PP _{<\lambda }$, the rest of the forcing being $(<\lambda
^+)$-closed, by the previous claim. (The power set of $\lambda$ is not
changed by $\PP_{\geq \lambda }\approx \PP_{>\lambda }$.)

\bigskip
\noindent
As is well-known, $j$ has the following properties:

\begin{description}
\item[ 1. ] crit($j$) = $\lambda$.
\item[ 2. ] $^\lambda M \subset M$.
\item[ 3. ] $\lambda^+ < j(\lambda ) < j(\lambda ^+) < \lambda ^{++}$.
\item[ 4. ] $M = \{j(F)(\lambda ) | F\in V  \wedge$ dom$(F) = \lambda \}$.
\end{description}

\noindent
We want to prove that $j$ `lifts' to an embedding
$\tilde j:V[G]\to N\subset V[G]$ (which we shall also denote by `$j$', abusing
notation), where $G$ is $\PP$-generic over $V$, thereby automatically
ensuring that $\lambda$ is also measurable in $V[G]$. 
The idea is to prolong the generic $G$ to a generic $H$ for the forcing
$j(\PP _{<\lambda })$.

\noindent
We start by comparing $\PP _{<\lambda }$ to $j(\PP _{<\lambda })$. $j(\PP
_{<\lambda })$ is an iteration defined in $M$, forcing the existence of
non-reflecting stationary sets on each $S^{\kappa^{++}}_\kappa$, for each
$\kappa$ successor inaccessible $<j(\lambda )$. If we compute the iteration
$j(\PP _{<\lambda })$ up to stage $\lambda$, we get exactly $\PP _{<\lambda }$.
We can now compute a generic extension $M[G]$ of $M$ using the $V$-generic
filter, since $\PP _{<\lambda }$-generics over $V$ are also generic over
$M$. Since $|\PP _{<\lambda }|<\lambda$, every canonical
$\PP _{\leq \lambda }$-name for a $\lambda$-sequence of ordinals
is in $M$, so in $V[G]$. Thus, we have

\medskip
\noindent
{\bf Claim:} $V[G]\models \mbox {  }^\lambda (M[G])\subset M[G]$.

\noindent
In $M[G]$, we will prolong $G$ to a $j(\PP _{<\lambda })$-generic. Call the
`remainder forcing' $\RR =\RR _{\lambda, j(\lambda )}$, that is,
$j(\PP _{<\lambda }) = \PP _{<\lambda }\ast \RR$. $\RR$ is
$(<\lambda _M^+)$-closed in $M[G]$; hence, by the previous claim, $\RR$ is
also $(<\lambda _M^+)$-closed in $V[G]$. 
In $M[G]$, the forcing $\RR$ is $j(\lambda )$-cc and has size $j(\lambda )$.
There are then at most $j(\lambda )$ maximal antichains in $\RR$ in the
model $M[G]$, as $[j(\lambda )^{j(\lambda )}=j(\lambda )]^M$, by
elementarity and since $|\PP _{<\lambda }|<\lambda$. In $V[G]$,
we can enumerate those antichains as
$\langle A_\alpha :\alpha <\lambda ^+\rangle$, since $\lambda^+ < j(\lambda ) <
j(\lambda ^+) < \lambda ^{++}$. Using the closure to meet all these
antichains, it is clear that in $V[G]$, we can build $H$ which is
$\RR$-generic over $M[G]$.
Letting $G^+ = G * H$, $G^+$ is $j(\PP _\lambda )$-generic over $M$. We define
$j:V[G]\to M[G^+]$ by $j(\dot\tau ^G) = j(\dot\tau)^{G^+}$. This is a
well-defined elementary embedding, as follows from the following fact,
whose proof uses the Truth Lemma and the elementarity of $j$.

\medskip
\noindent
{\bf Fact:} Let $k:M\to N$ be an elementary embedding between two transitive
models of ZFC. Let $\PP \in M$ be a forcing notion, let $k(\PP ) = \QQ$, and
suppose that $G$ is $\PP$-generic over $M$ and $H$ is $\QQ$-generic over
$N$. Let also $k``G\subset H$. Then the definition $k(\dot\tau ^G) =
k(\dot\tau)^H$ for every $\dot\tau \in M^\PP$ gives a well-defined elementary
embedding $k:M[G]\to N[H]$, which extends $k:M\to N$ and is such that
$k(G) = H$.

\medskip
\noindent
This ends the proof of measurability of $\lambda$ in $V[G]$, and thus the
proof of Theorem \ref{theorem:heightandgch}: At this point, to witness
the non-end-determining property of the appropriate $\kappa$'s, exactly the
same argument used in the proof of Theorem \ref{theorem:height} works.
Instead of using the failure/non-failure of $GCH$ at the successor
inaccessibles below $\kappa$ to code the corresponding object (club of
singulars or $\kappa$-Aronszajn tree), here, the coding is achieved
by the inclusion/non-inclusion of the corresponding $E$-set (see figure
\ref{picture:proof}).
\hfill $\Box_{\ref{theorem:heightandgch}}$

\bigskip

\section{Eees of Inner Models.}

\medskip
\noindent
In this final section, we push the study of the connections between
the existence of eees and the structure of the basic model one step
further: We look at inner models of some specific models related to $0^\sharp$
or to measures, and we look at some instances of the `wider model problem'.

\medskip
\noindent
In the last section, we remarked that for $L[\mu ]$, the situation is similar
to the one obtained via Theorem \ref{theorem:height}: if $\kappa$ is not
weakly compact, then there is a model $M$ of height $\kappa$ with no eees.
Contrasting this, we have the following fact for the inner model $L^M$, in
presence of $0^\sharp$.

\begin{proposition}\label{proposition:sharpsandL}
If $M\models `0^\sharp \mbox { exists'}$, and $M$ is
a {\it set\/} model of $ZFC$, then ${\cal E}_{L^M} \not= 0$.
\end{proposition}

\noindent
{\bf Proof:} The existence of
$0^\sharp$ in $M$ allows to code eees of $L^M$. \hfill $\Box$

\bigskip
\noindent
The converse to this is clearly false: just take $M = L_\kappa$ such that
(for example), $\kappa$ is weakly compact. Then, ${\cal E}_{L^M} \not= 0$,
yet $M\models `0^\sharp \mbox { does not exist'}$.

\medskip
\noindent
We have the following situation, asking the same question for {\it well
founded\/} eees (compare the following result and Proposition
\ref{proposition:sharpsandL}, and to Kunen's result in
Kaufmann [Ka 83]).

\begin{proposition}
Suppose that $M$ satisfies
\begin{itemize}
\item[ \bf i ]$M$ is countable, and
\item[ \bf ii ]$M\models `0^\sharp \mbox { exists'}$,
\end{itemize}

\noindent
and that $M$ has minimal ordinal height
among all the models which have simultaneously the properties {\bf i} and
{\bf ii}.
\underline{Then}, ${\cal E}_{L^M}^{wf}=0$.
\end{proposition}

\noindent
{\bf Proof:} We can use a $\Sigma_2^1$ formula to say that the height
of a certain model $M$ is $\alpha$ and that $M$ thinks that $0^\sharp$
exists: let $\psi (\alpha ) \equiv `\exists \mbox { transitive }M
(o(M)=\alpha \wedge M\models 0^\sharp \mbox { exists)'}$. This is clearly
$\Sigma_2^1$. Now, we also have by L\"owenheim-Skolem that
$ZFC \vdash [\exists \alpha \psi (\alpha )\to \exists \alpha <\omega _1\psi
(\alpha )]$.
The formula $\psi (\alpha )$ relativises down to models $N$ of $ZFC$
when $\alpha ^N$
is countable. So, working in $V$, we can fix $M$, a transitive model
such that $M\models
`0^\sharp$ exists', $o(M)=\delta$, and $\delta <\omega _1$ is the least
possible among the heights of such models. By Keisler-Morley, we know that
${\cal E}_{L^M}\not= 0$, and thus
$$L^M = L_\delta \models ZFC + \mbox { `inaccessibles are
cofinal in }ORD\mbox {'}.$$

\noindent
Suppose then that ${\cal E}_{L^M}^{wf}\not= 0$; let $\gamma <\omega _1$ be
such that $L_\gamma \succ L_\delta$. We then have that
$L_\gamma \models \aleph _\delta = \delta$. A collapse of $\delta$
to $\omega$ does the trick: working in $L_\gamma $, let $\PP = Coll(\omega
,\delta )$, and let $G$ be $\PP$-generic over $L_\gamma$.
So, on one hand, in $L_\gamma [G]$, $\delta$ is countable, and on the other
hand, $\psi (\delta )$ holds in $V$. Then, $\psi (\delta )$ also holds
in $L_\gamma [G]$, by $\Sigma _2^1$-absoluteness. But then, since
$L_\gamma$ and $L_\gamma [G]$ have the same ordinals, $\psi (\delta )$
holds in $L_\gamma$. Finally, by elementarity, $L_\delta \models
\exists \theta (\psi (\theta ))$. This contradicts the minimality
of $\delta$. \hfill $\Box$

\bigskip
\noindent
So far, the Height Problem has only been looked at for {\it inner\/} models
of a given ${\cal R}(\kappa )$ which is known to have non trivial eees. Yet
a natural question arises concerning `wider' models. More specifically, we
have the following

\begin{question} {\bf (eees for wider models)}
Let $\frak A$ be an inner model, ${\cal E}_{\cal R(\kappa )^{\frak A}}
\not= 0$,
$M\supset {\cal R}(\kappa )^{\frak A}$,
$M\models ZFC$, $o(M)=\kappa$, $|\kappa |>\omega$. When is ${\cal E}_M
\not= 0$?
\end{question}

\begin{figure}[ht]\hspace*{40.00mm}
\psfig{figure=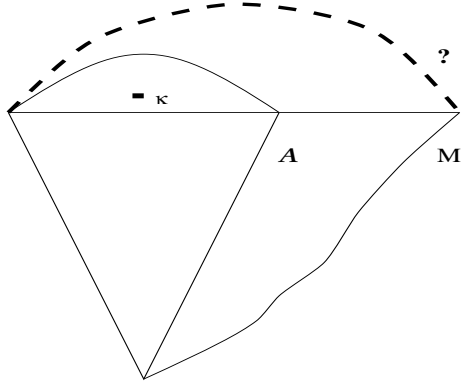,height=5cm,width=6cm}
\caption{ The `wider models problem' }
\end{figure}

\noindent
In such generality, there is no simple answer to this question. Yet,
we have that in particular, when $M$ is a generic extension of
${\cal R}(\kappa )^{\frak A}$, obtained via a {\it set-}forcing $\PP$, $M$
must also have
non trivial eees. This is easily established by observing that if $N\in
{\cal E}_{\cal R(\kappa )^{\frak A}}$, and $M={\cal R}(\kappa )^{\frak A}[G]$,
then $N[G]
\succ _eM$. This was implicitly used in the proof of Theorem
\ref{theorem:height}. There, we did not really have a set-forcing. Yet,
the decomposition properties of Easton forcing acting only on successor
inaccessible cardinals imply that in that case the forcing extension
necessarily has eees.

\noindent
An example of a CLASS forcing that destroys the property of ${\cal
R}(\kappa )$ having eees can be obtained from Boos's Easton style
forcing construction from~[Bo 74]. This
was pointed out to me by Ali Enayat. Boos's forcing, a variation on a
thin-set forcing
notion by Jensen, allows one to efface the Mahlo property:

\begin{theorem} (Boos - [Bo 74])
If $M$ is a transitive model of $ZFC + GCH$, and $\kappa$ is
$\kappa$-Mahlo in $M$, then for each $\alpha <\kappa$ there is a forcing
notion
$\BB _\alpha$ such that $M^{\BB _\alpha} \models `\kappa$ is
strongly
$\alpha$-Mahlo but not $(\alpha + 1)$-Mahlo'.
\end{theorem}

\noindent
Using the $\BB _1$ from this theorem, any weakly compact cardinal in $M$
becomes {\it the first\/} Mahlo in a generic extension. So, we have

\begin{corollary}
There are Easton-type forcings which do not
preserve the property of having eees: let $\kappa$ be a
weakly compact cardinal. Then,
although we know by Keisler and Silver [KeSi 70] that ${\cal E}_{{\cal
R}(\kappa )^M}\not= 0$, we have ${\cal E}_{{\cal R}(\kappa )^M[G]}= 0$,
for $G$ $\BB _1$-generic over $\cal R(\kappa )^M$.
\end{corollary}

\noindent
{\bf Remarks: 1)} Boos's construction ---a variant of a construction by
Jensen--- is an Easton-type iteration. By [KeMo 68], if $\kappa$ is the
first Mahlo cardinal, $\cal R(\kappa )$ cannot have any eees. The weakly
compact $\kappa$ of the corollary (in $M$) becomes the first Mahlo in the
extension by $\BB _1$.

\medskip
\noindent
{\bf 2)} Of course, if $\BB _n$ is used instead of $\BB _1$ here, then
$\kappa$ becomes the first $n$-Mahlo cardinal. The corresponding
${\cal R}(\kappa )^M[G]$ cannot either have
elementary end extensions: if $N\succ _e{\cal R}(\kappa )^M[G]$, then
$N\models \forall \lambda \exists C\subset \lambda (C$ club and $C$ does
not contain any $(n-1)$-Mahlos. But then, this is also true about $\kappa$,
and $\kappa$ would not be $n$-Mahlo.

\bigskip

\noindent
{\bf References.}
\begin{description}
\item[(Ba 76)] {\capit Baumgartner, J.}{ \it A new class of order types},
Annals of Mathematical Logic 9 (1976), 187-222.
\item[(BeJeWe 83)] {\capit Beller, A., Jensen, Welch}{ \it Coding the
Universe}, Cambridge Lecture Series, 1983.
\item[(Bo 74)] {\capit Boos, W.}{ \it Boolean Extensions which efface
the Mahlo property}, The Journal of Symbolic Logic, vol. 39, no 2, June 74,
pages 254-268.
\item[(CuD\v zSh 95)] {\capit Cummings, J., D\v zamonja, M., Shelah, S.}{
\it A consistency result on weak reflection}, accepted by Fund. Math.
\item[(En 84)] {\capit Enayat, A.}{ \it On certain Elementary Extensions of
Models of Set Theory}, Trans. Amer. Math. Soc., 1984.
\item[(En $\infty$)] {\capit Enayat, A.}{ \it Counting Countable
Models of Set Theory}, preprint.
\item[(Ka 83)] {\capit Kaufmann, M.}{ \it Blunt and Topless Extensions of
Models of Set Theory}, J. Symb. Logic, 48, 1983, 1053-1073.
\item[(KeMo 68)] {\capit Keisler, H.J., Morley, M.}{ \it Elementary
Extensions of Models
of Set Theory.} Israel J. of Math., vol. 6, 1968.
\item[(KeSi 70)] {\capit Keisler, H.J., Silver, J.} { \it
End Extensions of Models of Set Theory.}, Proc. Symp. Pure Math. 13
1970 (177-187).
\end{description}
\end{document}